%% file: main.tex
\documentclass{article}
\usepackage{import}
\usepackage{preambleCommands}

\title{Generalized Golub-Kahan bidiagonalization for nonsymmetric saddle point systems}

\author[1]{Andrei Dumitrasc}
\author[2]{Carola Kruse}
\author[1,2]{Ulrich R\"ude}

\affil[1]{Friedrich-Alexander-Universität Erlangen-Nürnberg, Erlangen, Germany}
\affil[2]{CERFACS, Toulouse, France}

\begin{document}

\maketitle

\begin{abstract}
    The generalized Golub-Kahan bidiagonalization has been used to solve saddle-point systems where the leading block is symmetric and positive definite. We extend this iterative method for the case where the symmetry condition no longer holds. We do so by relying on the known connection the algorithm has with the Conjugate Gradient method and following the line of reasoning that adapts the latter into the Full Orthogonalization Method. We propose appropriate stopping criteria based on the residual and an estimate of the energy norm for the error associated with the primal variable. Numerical comparison with GMRES highlights the advantages of our proposed strategy regarding its low memory requirements and the associated implications. 
\end{abstract}

\section{Introduction}

Solvers for systems of linear equations are at the core of many applications in scientific computing, often as part of a nonlinear solver. 
Direct solvers are typically employed for small, dense blocks, while iterative solvers deal with large, sparse matrices. Our following contribution is concerned with such iterative solvers, for the case of saddle point systems for which the (1,1)-block is nonsymmetric. The task of solving saddle point problems has a long history and many methods developed so far, for symmetric and nonsymmetric systems, have been collected in works such as \cite{benzi2005numerical,rozlovznik2018saddle}.

Saddle point systems arise in fields as diverse as incompressible fluid dynamics, structural mechanics, constrained and weighted least squares, and economics to name a few (see \cite{benzi2005numerical,pestana2015natural}). One famous example for a saddle point system coming from a partial differential equation (PDE), are the Navier-Stokes equations in computational fluid dynamics \cite{elman2014finite}. 
Depending on the application, 
the ensuing discretization leads to a large sparse system matrix.  
Such systems often have to be solved with  constraints regarding time or computing power, which motivates research in the field of efficient solvers.

Nonsymmetric problems are different to symmetric ones on a fundamental level, and solving them requires either significant modifications to existing symmetric solvers or completely new methods.
Sparse direct solvers can be employed to efficiently tackle nonsymmetric problems, but this strategy becomes infeasible when considering 3D industrial use cases, where the matrices involved are too large. The alternative here is to use iterative solvers, such as Krylov subspace solvers combined with various preconditioners (see \cite{Farahbakhsh2020_KrylovSubspaceMethods,Meurant2020_KrylovMethodsNonsymmetric,ordonez2023scalable,scott2014signed}). Iterative methods are attractive 
due to their customizable character: 
one has a finer degree of control over the complexity, requirements and quality of the solution when compared to a direct solver. 
For instance, when the linear iterative solver is a building block of a nonlinear one, the solution the former has to deliver typically needs only a few accurate digits, meaning the iterative solver can be stopped early. While iterative methods require significantly less memory than direct ones, 3D nonsymmetric problems can still prove challenging. Iterative solvers from the Krylov family approach the issue of nonsymmetry in two different ways, leading to two groups of algorithms (see \cite{saad2003iterative,barrett1994templates}). One of them, with \ac{BiCG} as a concrete example, uses relatively inexpensive iterations with a constant cost, but irregular convergence. The other group, with \ac{GMRES} as example, involves iterations which get more and more expensive, as a result of storing and processing more and more vectors. The resulting advantage is then smooth, predictable convergence behavior. Restarts can be integrated into such algorithms to limit memory usage, but this also lead to slower convergence. We see then how all these solution approaches represent different trade-offs between  a solver's computational requirements and its speed or accuracy.

Solvers for the kind of indefinite saddle point system, may or may not exploit the 2-by-2 block structure of the considered matrix. The solver can either be applied in an all-at-once manner, i.e., to the system as a whole, or in a segregated way, where the operations use the block structure. 
All-at-once solvers, that can be applied to such indefinite problems are, e.g., sparse direct methods such as MUMPS \cite{MUMPS01}, preconditioned Krylov subspace solvers \cite{saad2003iterative} (for example \ac{GMRES}, \ac{BiCG}, CGS, CGNE)  or also a monolithic multigrid as solver or preconditioner (\cite{awanou2005convergence}, see \cite{kohl2022textbook} for symmetric Stokes flow). Examples for segregated algorithms are stationary methods, Uzawa-like methods \cite{uzawa58} or the Augmented-Lagrangian-Uzawa algorithm \cite{fortin83augmentedlagrangian}.
The solver we propose in this work also falls into the class of segregated methods. We build our development on a generalized Golub-Kahan bidiagonalization (introduced by Mario Arioli in \cite{arioli2013generalized}) which has been successfully used for symmetric saddle point problems (see \cite{KrSoArTaRu2021,KrDaTaArRu2020,KrSoArTaRu2020}). However, the potential of this algorithm has not yet been developed and exploited in the nonsymmetric case.

In our previous work (see \cite{KrSoArTaRu2021,KrDaTaArRu2020,KrSoArTaRu2020}), we have tackled symmetric systems using the generalized CRAIG solver \cite{arioli2013generalized}. In this paper, our aim is to extend this algorithm to the nonsymmetric case. The setting for an equivalence between CRAIG and CG has been described in \cite[Chapter~5]{orban2017iterative}, and adapting it to include nonsymmetry is the first step towards deriving the algorithm we propose. The core theoretical component is an equivalence with the Arnoldi process applied to the Schur complement, which we describe in detail in \Cref{sec:schurComp}. 
The associated solver is then equivalent to FOM, in a similar way to the \ac{CG} equivalence explained in \cite[Chapter~5]{orban2017iterative}. 

Given the iterative nature of our solver, it is important to define stopping criteria, especially if we are only interested in a solution of moderate accuracy. 
To that end, we give inexpensive means of computing the residual norm, 
as well as an estimate of the error norm. 
We highlight the practical advantages of our solver by a number of comparisons in a numerical setting.
We generate Navier-Stokes problems using Incompressible Flow \& Iterative Solver Software \footnote{http://www.cs.umd.edu/~elman/ifiss3.6/index.html} (IFISS).  The resulting nonlinear systems are solved using Picard's iteration, which relies on a linear solver to compute its correction. 
In this context, we compare our solver with \ac{GMRES} and FOM.

Our paper is organized as follows. In \Cref{sec:schurComp} we analyze the existing setting for symmetric saddle point problems. In \Cref{sec:nsGKB} we extend the considerations in the previous section to the nonsymmetric case and introduce an adapted Golub-Kahan bidiagonalization, the corresponding linear solver and its stopping criteria. Next, we numerically compare our solver against other choices in \Cref{sec:tests}. Finally, we conclude with summarizing remarks and future development ideas in \Cref{sec:conc}. \\

\section{Connections with the Schur complement}
\label{sec:schurComp}
In this article we will study the iterative solution of saddle point systems of the form
\begin{equation}\label{eq:ini_sys}
	\left[
	\begin{matrix}
		\Mm & \Am \\
		\Am ^T & \mZ 
	\end{matrix}
	\right]
	\left[
	\begin{matrix}
		\uv  \\
		\pv 
	\end{matrix}
	\right]
	=
	\left[
	\begin{matrix}
		\mZ  \\
		\bv 
	\end{matrix}
	\right], 
\end{equation}
where $ \Mm \in \bR ^{m \times m}, \Am \in \bR ^{m \times n} $ and  $ \bv \in \bR ^n$. 
We assume that $\Mm$ is nonsymmetric ($\Mm \neq \Mm ^T$) and positive definite ($ \xv^T \Mm  \xv > 0, \forall \, \xv \in \bR ^m $) and that $\Am$ has full column rank, following the respective definitions for real matrices given in i.e. \cite{saad2003iterative}. 
If $\Mm$ is positive semidefinite and $\mbox{ker}(\Mm) \, \cap \,  \mbox{ker}(\Am^T)=\{ \mZ \}$, we can use the augmented Lagrangian method as in \cite{KrDaTaArRu2020,KrSoArTaRu2021}, which leads to a positive definite matrix.
The more general case
\begin{equation}\label{eq:gen_sys}
	\left[
	\begin{matrix}
		\Mm & \Am \\
		\Am ^T & \mZ  
	\end{matrix}
	\right]
	\left[
	\begin{matrix}
		\wv  \\
		\pv 
	\end{matrix}
	\right]=
	\left[
	\begin{matrix}
		\bv _1\\
		\bv _2 
	\end{matrix}
	\right] 
\end{equation}
can be transformed into \eqref{eq:ini_sys} by first setting $ \wv _0 = \Mm ^{-1} \bv _1$ and then considering the residual equation corresponding to the point $[\wv _0 \ \mZ ] ^T$. This leads to the form \eqref{eq:ini_sys}, where $ \uv = \wv - \wv _0 $ and $ \bv = \bv _2 - \Am ^T  \wv _0  $.

The assumption about $\Mm$ being positive definite is indeed a realistic, as among applications leading to nonsymmetric saddle point matrices, there is a subgroup for which this is true. The presence or absence of this property stems from the interplay of physical parameters, linearization, discretization and possibly others. One such example is the Navier-Stokes problem linearized with Picard’s iteration (Oseen problem), and discretized with the Compatible Discrete Operator (see \cite{jang2023fast}) or with Taylor-Hood finite elements. In \cite{konshin2015ilu}, a discussion is presented on the choice of certain parameters to obtain a positive definite nonsymmetric (1,1)-block. An alternative if this condition does not hold (for symmetric systems) is the null-space approach, such as that explored in \cite{scott2022null} for saddle point systems with rank deficient blocks.

One way of solving \Cref{eq:ini_sys} is via the Schur complement equation. 
We can find $\pv$ as $\pv = - \Sm ^{-1} \bv,$ with $\Sm = \Am ^T  \Mm ^{-1} \Am $, by substituting $\uv$ in the second equation, based on the first. 
Then, we have $\uv = - \Mm ^{-1} \Am \pv$.

We base our following development on a generalized Craig algorithm for symmetric saddle point matrices described by M. Arioli in \cite{arioli2013generalized} and assume for the rest of this section symmetry of $\Mm$. 
In \cite[Chapter~5]{orban2017iterative}, the following important connection has been established: applying the generalized \ac{GKB} to a saddle-point system with a matrix having a symmetric leading block is equivalent to applying \ac{CG} to the associated Schur complement equation. 
As such, we can draw a link between $\Am$'s bidiagonalization and $\Sm$'s tridiagonalization given by the Lanczos process specific to CG.

We consider the generalized Golub-Kahan bidiagonalization of $\Am$ in the form
\begin{equation}
\label{eq:symGKB}
    \begin{cases}
    \Am \Qm  = \Mm \Vm  \Bm , &\qquad \Vm  ^T \Mm \Vm   = \Id,  \\
    \Am ^T \Vm  = \Qm  \Bm  ^T ,  &\qquad \Qm  ^T \Qm  = \Id, 
    \end{cases}
\end{equation} 
with the upper bidiagonal $ \Bm \in \bR ^{n \times n}$, the orthogonal $ \Qm \in \bR ^{n \times n}$ and $ \Vm \in \bR ^{m \times n}$.
The basis $\Vm$ has $\Mm$-orthogonal columns with respect to the inner product and norm
\begin{equation*}
 \innprodM{\xv}{\yv}= \xv ^T \Mm \yv, \qquad 
    \normM{\xv}  = \sqrt{ \innprodM{\xv}{\xv} }.
\end{equation*}
With a starting vector $\bv \in \bR  ^n$, the algorithm first computes $\beta _1 = \normEu{\bv} $ and sets $\qv _1 = \bv / \beta _1$, which is the first element of the orthogonal basis $\Qm$. 
The partial decomposition at step $k$ is given iteratively by
\begin{equation}
\label{eq:sGKB}
    \begin{cases}
    \Am \Qm _k = \Mm \Vm _k \Bm _k, &\qquad \Vm _k ^T \Mm \Vm _k  = \Lm _k,   \\
    \Am ^T \Vm _k = \Qm _k \Bm _k + \beta _ {k+1} \qv _{k+1} \ev _k^T,  &\qquad \Qm _k ^T \Qm _k = \Id _k,
    \end{cases}
\end{equation}
with
\begin{equation}
\label{eq:Bsymm}
    \Bm _k=
    \left[ 
	\begin{matrix}
	\alpha_1 & \beta_2 & 0 & \ldots & 0 \\
    0 &	\alpha_2 & \beta_3 &  \ldots & 0 \\
    \vdots &	\vdots & \vdots & \vdots & \vdots  \\
    0 & \ldots & 0 & \alpha_{k-1} & \beta_k \\
    0 & \ldots & 0 & 0 & \alpha_{k}
	\end{matrix}
	\right],
 \end{equation}

Reconsidering the equivalence of CRAIG with \ac{CG} for the Schur complement, from \Cref{eq:symGKB} we find the tridiagonalization of $\Sm$
\begin{align*}
    \Sm &= \Am ^T  \Mm ^{-1} \Am \\
        &= ( \Mm \Vm  \Bm \Qm ^T ) ^T \Mm ^{-1}( \Mm \Vm  \Bm \Qm ^T ) \\
        &= \Qm ( \Bm ^T \Bm ) \Qm ^T,
\end{align*}
where $( \Bm ^T \Bm )$ is the tridiagonal matrix. 
The relationship between these two processes also holds when considering only partial bidiagonalization and tridiagonalization as in \Cref{eq:sGKB}, i.e., when one performs less than $n$ iterations.

Using this 
equivalence, we can implicitly tridiagonalize the Schur complement in the sense that we find $\Qm$ and $\Bm $ without knowledge of $\Sm$. 
It is important to note that the Lanczos process leads to a tridiagonal matrix only if the target matrix (in this case, $\Sm$) is symmetric. 
This condition opens the possibility of short recurrences, where one vector needs to be orthogonalized only against the previous two. 
Then, in exact arithmetic, it will also be orthogonal to all the previous ones. 
Consequently, only the most recent vectors are stored and used.

\section{The \ac{GKB} for nonsymmetric systems}
\label{sec:nsGKB}

We turn to the case where $\Mm$ is not symmetric. 
It follows that the resulting Schur complement $\Sm$ is also nonsymmetric and we can no longer use \ac{CG} to solve such a system, leading us to consider the \ac{FOM}, see \cite{saad2003iterative}.
%

In our case, 
the following discussion will be based on FOM, 
which is motivated by the fact that \ac{FOM} is more similar to \ac{CG}, which gives the equivalence to the CRAIG algorithm in the symmetric case. Our assumption that $\Mm$, the leading block in \Cref{eq:ini_sys}, is positive definite is important for the convergence of \ac{FOM} since it ensures that no breakdown can occur. This point is discussed further in \cite{frommer2020krylov}, with an emphasis on the relationship between the (quadratic) numerical range of a matrix and the way FOM and \ac{GMRES} converge. The other reason why we require positive definiteness for $\Mm$ is inspired by the symmetric version of the generalized CRAIG, where the left vectors have unit $\Mm$-norm.

Choices other than \ac{FOM} could be algorithms that deal with nonsymmetry using biorthogonalization, such as \ac{BiCG} and its successors (see \cite{saad2003iterative,barrett1994templates}). For example, when solving nonlinear equations, one can employ a linearization via Chebyshev interpolation, with the resulting system to be solved by \ac{BiCG}, as done in \cite{correnty2022preconditioned} for parameterized systems. If using \ac{GMRES}, one can construct preconditioners as in \cite{loghin2006bounds}, or more recent ones like \cite{cao2015relaxed,wei2022new}, based on a Hermitian and skew-Hermitian splitting.

In \ac{FOM}, the Arnoldi process generates mutually orthogonal vectors using $\Sm$, but 
with long recurrences. 
In this case, we need to store all the vectors and generate a new one by orthogonalizing against all of them.
Therefore, instead of achieving a tridiagonalization of $\Sm$, as in the symmetric case, we only reduce it to an upper Hessenberg form. 
We still want to find this in an implicit way, without referring to $\Sm$, but rather to $\Am$. 
Then, it is necessary to redefine the Golub-Kahan process as follows. 

With a starting vector $\bv \in \bR  ^n$, we first set $\qv _1 = \bv / \beta _1$ with $\beta _1 = \normEu{\bv} $, which is the first element of the orthogonal basis $\Qm$.
The partial decomposition at step $k$ is given recursively by
\begin{equation}
\label{eq:nsGKB}
    \begin{cases}
    \Am \Qm _k = \Mm \Vm _k \Bm _k, &\qquad \Vm _k ^T \Mm \Vm _k  = \Lm _k,   \\
    \Am ^T \Vm _k = \Qm _k \Hm _k + \beta _ {k+1} \qv _{k+1} \ev _k^T,  &\qquad \Qm _k ^T \Qm _k = \Id _k,
    \end{cases}
\end{equation}
with $\Bm_k \in \bR ^{k \times k}$ given in \Cref{eq:Bsymm} and 
\begin{equation}
\label{eq:BandHmat}
	\Hm _k=
    \left[ 
	\begin{matrix}
\theta_{1,2}	&	\theta_{1,3}	&	\theta_{1,4}	&	\ldots	&	\theta_{1,k+1}	&	\\
\beta_2	&	\theta_{2,3}	&	\theta_{2,4}	&	\ldots	&	\theta_{2,k+1}	&	\\
0	&	\beta_3	&	\theta_{3,4}	&	\ldots	&	\theta_{3,k+1}	&	\\
\vdots	&	\vdots	&	\vdots	&	\vdots	&	\vdots	&	\\
0	&	0	&	0	&	\beta_{k}	&	\theta_{k,k+1}	&	
	\end{matrix}
	\right],
\end{equation}
where $ \Qm _k \in \bR ^{n \times k}, \Vm _k \in \bR ^{m \times k}, \Hm _k, \Lm _k$ and $ \Id _k \in \bR ^{k \times k}$, the latter being the identity matrix. The vector $\ev _k$ is the last column of $\Id _k$. In contrast to \Cref{eq:sGKB}, we now have the Hessenberg matrix in the second equation, whereas the first equation stays unchanged.

\begin{proposition}
\label{prop:KlowTri}
The matrix $\Lm _k$ is unit lower triangular.
Also, $ \Bm _k^T \Lm ^T_k= \Hm_k$. As such, we can think of  $ \Bm _k^T $ and $\Lm _k^T$ as the LU factors of $\Hm_k$, respectively.
\end{proposition}
\begin{remark}
\label{prop:QorthoRes}
The vector $\qv _{k+1}$ is defined such that $\Qm _k ^T \qv _{k+1} = 0 $. 
In exact arithmetic, the decomposition finishes at the latest at step $k=n$, when the term $\beta _ {k+1} \qv _{k+1} \ev _k^T$ vanishes. 
\end{remark}

Given an $n$-step decomposition of type \eqref{eq:nsGKB}, in light of \Cref{prop:KlowTri} and \Cref{prop:QorthoRes}, we can find the reduction to upper Hessenberg form of the Schur complement matrix $\Sm$ as 
\begin{equation}\label{eq:GKBequivSchur}
    \begin{split}
    \Sm &= \Am ^T  \Mm ^{-1} \Am \\
        &= ( \Mm \Vm  \Bm \Qm ^T ) ^T \Mm ^{-1}( \Mm \Vm  \Bm \Qm ^T ) \\
        &= \Qm ( \Bm ^T \Lm ^T \Bm ) \Qm ^T \\
        &= \Qm ( \Hm \Bm ) \Qm ^T \\
        &= \Qm  \Hm _{\Sm}  \Qm ^T ,
    \end{split}
\end{equation}
with $\Hm _{\Sm} = \Hm \Bm $ being an upper Hessenberg matrix. As for the symmetric case, the connection between the two processes is still valid for $k<n$. The solution $\pv$ is found as in \ac{FOM} \begin{equation}\label{eq:solSchurp}
    \pv = - \Qm \Hm _{\Sm} ^{-1} (\beta _1 \ev _1),
\end{equation}
then replaced in the first row of \eqref{eq:ini_sys} to yield
\begin{equation}\label{eq:solSchuru}
    \uv = - \Mm ^{-1} \Am \pv.
\end{equation}
The presence of the minus sign in \eqref{eq:solSchurp} is a consequence of the choice of starting point: the \ac{GKB} starts with $\bv$, while \ac{FOM} on the Schur complement starts with $-\bv$.

\subsection{The nonsymmetric CRAIG algorithm}
Let us now define the nonsymmetric GKB decomposition, which will be used to define our proposed nonsymmetric solver, nsCRAIG. The matrices in \Cref{eq:nsGKB} are built step by step, as follows. With an initial vector $\qv _1 = \bv / \beta _1$ with $\beta _1 = \normEu{\bv}$, we can compute the matrices $\Vm$ and $\Qm$ step by step by using the relation in \Cref{eq:nsGKB}. 
The first left vector, $\vvv _1$, is found by
\begin{equation*}
    \begin{split}
    \wv _1 &= \Mm ^{-1} \Am \qv _1, \\
    \alpha _1 &= \normM{\wv _1}, \\
    \vvv _1 &= \wv _1 / \alpha _1 .
    \end{split}
\end{equation*}
For the normalization of the left vectors in general, it is necessary to have $\wv ^T \Mm \wv > 0, \forall \ \wv \in \bR ^m, \wv \neq \mZ$, which follows from our requirement that $\Mm$ is positive definite.

For the subsequent steps, a new right vector is computed following
\begin{equation*}
    \begin{split}
    \hat{\gvv} _k &= \Am ^T \vvv _{k-1}, \\
    \hv _k &= \Qm _{k-1} ^T \hat{\gvv} _k,\\
    \gvv _k &= \hat{\gvv} _k - \Qm _{k-1} \hv _k .
    \end{split}
\end{equation*}
This way, we orthogonalize the $k$-th right vector against all the previous ones stored in $\Qm _{k-1}$.
For simplicity, we have expressed this using the Gram-Schmidt method. 
However, the vectors can rapidly lose orthogonality in practical applications. 
For improved numerical reliability, we can replace the Gram-Schmidt method with Modified Gram-Schmidt, orthogonal transformations such as Givens rotations  \cite{saad2003iterative} or newer methods involving random sketching \cite{balabanov2022randomized}.
The vector $ \hv _k = [\theta _{1,k}, \theta _{2,k}, ..., \theta _{k-1,k} ] ^T$ is the $k-$th column of $\Hm$, consisting of the entries on the diagonal and above. 
This vector stores the inner products 
that are needed to orthogonalize $\hat{\gvv} _k $ against
all the 
previous 
vectors in $\Qm$. 
The resulting right vector is then normalized and stored as
\begin{equation*}
    \begin{split}
    \beta _k &= \normEu{\gvv _k} \\
    \qv _k &=  \gvv _k / \beta _k \\
    \Qm _k &= [\Qm _{k-1} \ \qv _k ].
    \end{split}
\end{equation*}
Note that it is necessary to store all the right vectors and use them in the orthogonalization process, to maintain global mutual orthogonality.

The left vectors are computed as
\begin{equation*}
    \begin{split}
    \wv _k &= \Mm ^{-1} ( \Am \qv _k - \beta _k \Mm  \vvv _{k-1})\\
    \alpha _k &= \normM{\wv _k} \\
    \vvv _k &= \wv _k / \alpha _k .
    \end{split}
\end{equation*}
%
%
Here, it is sufficient to store only the latest left vector $\vvv _{k-1}$ and use it to compute $\vvv _k$
.
Thus the memory cost of the algorithm is the storage of the full basis of the right-hand side vectors of size $n$ and a three-term recurrence for the left-hand side vectors of size $m$.

This together results in the nonsymmetric CRAIG algorithm, nsCRAIG, as presented in \Cref{alg:nsGKB1}.
\begin{algorithm}[H]
  \caption{nsCRAIG algorithm.}
  \label{alg:nsGKB1}
  \begin{algorithmic}
  \REQUIRE \hspace{-0.2cm}: {$\Mm , \Am , \bv$, maxit, tol}
  \ENSURE \hspace{-0.2cm}: $\uv _k, \pv _k$
  \vspace{.5cm}
  
  \STATE{$\beta_1 = \normEu{\bv}; \qquad$ $\qv_1 = \bv / \beta _1;$}
  \STATE{$\wv = \Mm^{-1} \Am \qv_1 ;\qquad $ $\alpha_1 = \normM{\wv} ;\qquad $ $\vvv = \wv / \alpha_1;$}
  \STATE{$\chi _1 = \beta _1 / \alpha _1 ;\qquad $ $ \rho = 1 / \alpha _1 ;\qquad $ $ k = 1 ;$ }  
  \WHILE{ $k \leq $ maxit}
  \STATE{$\gvv = \Am^T \vvv; \qquad$ $\hv _k = \Qm ^T \gvv; \qquad$  $\gvv = \gvv - \Qm \hv _k; $  }
  \STATE{$\beta_{k+1} = \normEu{\gvv}; \qquad $    }
  \IF{ \textit{Stopping criterion}$(\chi_k,\rho)$ = TRUE} 
  \STATE{ \textbf{break;} }
  \ENDIF
  \STATE{$\qv_{k+1} = \gvv / {\beta_{k+1}}; \qquad $ }
  \STATE{$\wv = \Mm ^{-1} \left(  \Am \qv _{k+1} - \beta _{k+1} \Mm \vvv \right);\qquad $ }
  \STATE{$\alpha _{k+1} = \normM{\wv};\qquad $ $\vvv  = \wv / {\alpha _{k+1} }$}
  \STATE{$\chi _{k+1} = - \dfrac{\beta _{k+1}}{\alpha _{k+1}} \chi _k; \qquad $ $k = k + 1;$}
  \ENDWHILE
  
  \STATE{$ \Hm _{\Sm} = \Hm \Bm ; \qquad $  $ \yv  = \Hm _{\Sm} ^{-1}  ( \beta _1  \ev _1 ) ; \qquad $ $ \pv _k = \Qm \yv  ; \qquad $  $  \uv _k = - \Mm ^{-1} \Am \pv _k; $ }
  \end{algorithmic}
\end{algorithm}
Although not yet introduced, we mention here that computing the quantities $\rho$ and $\chi$ is necessary for the stopping criterion, 
which we will explain in \Cref{subsec:stopMns}. 

\subsection{Linear system solution}

So far, we have 
described how to build $\Vm$ and  $\Qm$, the left- and right-hand bases, if the (1,1)-block $\Mm$ is nonsymmetric. 
Our primary interest is, however, to solve the linear system in \Cref{eq:gen_sys}, which we will describe next and thus explain the last line in \Cref{alg:nsGKB1}. 
We assume exact arithmetic and that we have a complete (i.e. $k=n$)  decomposition of type \eqref{eq:nsGKB}. 
For simplicity, we omit the subscript for quantities coming directly from this decomposition.

\begin{lemma}\label{lem:nsGKBsol}
The solution $ [\uv \ \pv] ^T $ to \Cref{eq:ini_sys} can be computed as
\begin{equation*}
    \begin{cases}
    \uv = \Vm \zv, \qquad \zv = \beta _1 \Lm ^{-T} \Bm ^{-T} \ev _1, \\
    \pv = \Qm \yv, \qquad \yv = -\Bm ^{-1} \zv .
    \end{cases}
\end{equation*}
\end{lemma}
\begin{proof}
The first step is to extend $ \Vm$ with $m-n$ columns, such that 
\begin{equation}
\label{eq:nonsym_ortho}
    \begin{cases}
    \innprodM{\vvv _i}{\vvv _j} = 0, & \qquad \forall i \in [1,m-1], \ \forall j \in [i+1,m],  \\
    \normM{\vvv _i} = 1, & \qquad \forall i \in [1,m].
    \end{cases}
\end{equation}
This condition is inspired by the case where $\Mm$ is symmetric and positive definite, and can be used to define an inner product and norm. With these, one can define an $\Mm-$orthogonal basis. 
For our case with nonsymmetric but still positive definite $\Mm$, this is not possible. Nonetheless, condition (\ref{eq:nonsym_ortho}) is sufficient for the purposes of this proof, and can be intuitively considered a kind of one-sided, nonsymmetric orthogonality.

Let us define $\Vm _1 = \Vm $. By this, we mean that $\Vm$ represents the first part of the following matrix $\bar{ \Vm } = [\Vm _1 \ \Vm _2  ] \in \bR ^{m \times m}$. The second part, $\Vm _2$, is a set of additional $m-n$ columns, formally appended after the process in \eqref{eq:nsGKB}, such that the extended matrix is of order $m$. 

Additionally, we define $ \Lm _m = \bar{ \Vm } ^T \Mm \bar{ \Vm }  $, where
\begin{equation*}
    \Lm _m =
    	\left[
	\begin{matrix}
		 \Lm _{1,1} & \mZ _{n,m-n} \\
		\Lm _{2,1} & \Lm _{2,2} 
	\end{matrix}
	\right],
	 \Lm _m \in \bR ^{m \times m},
\end{equation*}
with the individual blocks
\begin{align*}
    \Lm _{1,1} &=  \Vm _1 ^T \Mm \Vm _1, \qquad   \Lm _{1,1} \in \bR ^{n \times n}, \text{ unit lower triangular,}   \\
    \Lm _{2,1} &=  \Vm _2 ^T \Mm \Vm _1, \qquad   \Lm _{2,1} \in \bR ^{m-n \times n}   \\
    \Lm _{2,2} &=  \Vm _2 ^T \Mm \Vm _2, \qquad   \Lm _{2,2} \in \bR ^{m-n \times m-n}, \text{ unit lower triangular} .   
\end{align*}
The block $ \Vm _2 $ is added after the bidiagonalization in \eqref{eq:nsGKB}, so the other matrices therein are not changed, merely extended as
\begin{equation}\label{eq:ext_nsGKB}
    \begin{cases}
    \Am \Qm = \Mm \bar{\Vm} 
        	\left[
	\begin{matrix}
		 \Bm \\
		\mZ _{m-n,n}
	\end{matrix}
	\right],
    \\
    \Am ^T \bar{\Vm}  = \Qm  [\Hm \ \mZ _{n,m-n} ].
    \end{cases}
\end{equation}
If we left-multiply the first row with $\bar{\Vm}$ and the second row with $\Qm ^T$, we get \begin{equation*}
        	\left[
	\begin{matrix}
		 \Lm _{1,1} & \mZ _{n,m-n} \\
		\Lm _{2,1} & \Lm _{2,2} 
	\end{matrix}
	\right]
	\left[
	\begin{matrix}
		 \Bm \\
		\mZ _{m-n,n}
	\end{matrix}
	\right]
	=
	\left[
	\begin{matrix}
		 \Hm ^T \\
		\mZ _{m-n,n}
	\end{matrix}
	\right],
\end{equation*}
which can be seen as a version of \Cref{prop:KlowTri} that uses the extended matrices we defined above. 
 It is of course still true that  $\Lm _{1,1} \Bm = \Hm ^T  $, as before the extension with $ \Vm _2 $. Additionally, from the second row, we get 
 \begin{equation}
     \label{eq:LB0} 
 \Lm _{2,1} \Bm = \mZ _{m-n,n} .
 \end{equation}
%
%
We next introduce a basis transformation for the system \eqref{eq:ini_sys}, using a block-diagonal matrix with diagonal elements $ \bar{ \Vm }  $ and $\Qm $. The resulting system is 
\begin{equation}\label{eq:transfSys}
	\left[
	\begin{matrix}
		 \Lm _{1,1} & \mZ _{n,m-n} &  \Hm  ^T \\
		\Lm _{2,1} & \Lm _{2,2} & \mZ _{m-n,n} \\
		\Hm  & \mZ _{n,m-n} & \mZ _{n,n}
	\end{matrix}
	\right]
	\left[
	\begin{matrix}
		\zv _1  \\
		\zv _2  \\
		\yv 
	\end{matrix}
	\right]
	=
	\left[
	\begin{matrix}
		\mZ _{n,1} \\
		\mZ _{m-n,1} \\
		\Qm  ^T \bv 
	\end{matrix}
	\right], 
\end{equation}
which follows by examining the individual blocks and their relationships based on \eqref{eq:ext_nsGKB}.
The solution to the initial system \eqref{eq:ini_sys} can then be obtained as 
\begin{equation*}
    \left[
	\begin{matrix}
		\uv \\
		\pv 
	\end{matrix}
	\right]
	=
    \left[
	\begin{matrix}
		\bar{ \Vm }  & \mZ _{m,n}  \\
		\mZ _{n,m} & \Qm 
	\end{matrix}
	\right]
	\left[
	\begin{matrix}
		\bar{ \zv } \\
		\yv 
	\end{matrix}
	\right]
	=
		\left[
	\begin{matrix}
		\Vm _1 \zv _1 +  \Vm _2 \zv _2 \\
		\Qm \yv 
	\end{matrix}
	\right]
	, 
\end{equation*}
with the additional notation $ \bar{ \zv } = [ \zv _1 \  \zv _2 ] ^T $.
We proceed to identifying $\zv _1, \zv _2$ and $\yv $. 
From the third row of \eqref{eq:transfSys}, we find 
\begin{equation*}
    \zv _1 = \Hm  ^{-1} \Qm  ^T \bv.
\end{equation*}
%
%

With \Cref{prop:KlowTri} and noting that the first step of the decomposition in \eqref{eq:nsGKB} is to set $\beta _1 =  \normEu{\bv}$ and $\qv _1 = \bv / \beta _1$, we may rewrite
\begin{equation*}
    \zv _1 = \beta _1 \Lm _{1,1} ^{-T} \Bm  ^{-T}  \ev _1.
\end{equation*}
Next, we obtain from the first row of \eqref{eq:transfSys} and \Cref{prop:KlowTri}
\begin{equation*}
    \yv  = -  \Hm  ^{-T} \Lm _{1,1} \zv _1=-  \Bm  ^{-1} \zv _1.
\end{equation*}
From the second row of \eqref{eq:transfSys} and the previous equation follows
\begin{equation*}
    \zv _2 = - \Lm _{2,2} ^{-1} \Lm _{2,1}  \zv _1=  \Lm _{2,2} ^{-1} \Lm _{2,1} \Bm    \yv.
\end{equation*}
Since $ \Lm _{2,1} \Bm  = \mZ _{m-n,n} $, 
it follows that $ \zv _2 =  \mZ _{m-n,1}$. Hence, the corresponding vectors in $\Vm _2$, introduced by the extension defined at the beginning of the proof, are not necessary. This means that we can recover $\uv$ only using $\Vm _1$, which are the $n$ left vectors provided by the decomposition in \eqref{eq:nsGKB} and $ \zv _1 $. 
Finally, since only $\zv _1$ and  $\Vm _1$ are relevant, we drop the subscript $1$. 
\end{proof}

We conclude with the following remarks. 
For exact arithmetic, the nsCRAIG algorithm gives the exact solution after at most $n$ steps. 
In the absence of exact arithmetic and/or if we only have a partial decomposition ($k<n$), the formulas in \Cref{lem:nsGKBsol} give only an approximate solution. 
The bidiagonalization progresses by relying on the vectors $\vvv $ and $ \qv$, without referring to the current iterates $\uv _k$ and $ \pv _k$.
Consequently, we can postpone forming this approximation until we have found that the stopping condition has been fulfilled. 
Then, we only need to apply the inverses $\Lm _k ^{-T}$ and $\Bm _k ^{-1}$  mentioned in \Cref{lem:nsGKBsol} once. 
Owing to the particular structure of the matrices  $\Bm _k ^T$ and $\Lm _k ^T$, their inverses can be applied by using forward and backward substitution.

\subsection{Stopping criteria}
\label{subsec:stopMns}

In practical applications, Krylov solvers do not aim to achieve an exact solution, but rather a sufficiently accurate one. The solver stops before performing $n$ iterations, by monitoring a quantity which is indicative of the accuracy of the current iterate. One common choice for this quantity is the residual, more precisely its Euclidean norm. After $k$ steps of the decomposition in \eqref{eq:nsGKB}, we have an approximation $\uv _k = \Vm _k \zv _k$, with a corresponding residual  $\| \bv - \Am ^T \uv _k \| _2$. As we can see from \Cref{lem:nsGKBsol}, computing $\uv _k$ involves expenses such as inverting matrices and storing all the vectors from $\Vm$. With this in mind, we would like a simpler expression for the norm, independent of $\uv _k$. As for other Krylov solvers \cite{saad2003iterative}, this is possible and can be achieved with just a few scalar operations, as follows. 

\begin{lemma}
\label{lem:resNorm}
    The residual norm at step $k$  is  
 \begin{equation*}
    \| \bv - \Am ^T \uv _k \| _2= \beta _ {k+1} \chi _k,
\end{equation*}   
where 
\begin{equation*}
    \chi _1 =  \frac{\beta _1}{\alpha _1},  \qquad
    \chi _k = - \frac{\beta _k}{\alpha _k}  \chi _{k-1}.
\end{equation*}
\end{lemma}
\begin{proof}
Using \Cref{lem:nsGKBsol} and \Cref{eq:nsGKB}, we find
\begin{equation*}
    \begin{split}
    \bv - \Am ^T \uv _k 
     &= \beta _1 \Qm _k \ev _1 - \Qm _k \Hm _k \zv _k - \beta _ {k+1} \qv _{k+1} \ev _k ^T \zv _k
    \end{split}
\end{equation*}
and  
\begin{equation*}
    \Qm _k \Hm _k \zv _k = \beta _1 \Qm _k \Hm _k \Hm _k ^{-1} \ev _1 = \beta _1 \Qm _k \ev _1.
\end{equation*}
It follows
\begin{equation*}
     \bv - \Am ^T \uv _k = -\beta _{k+1} \qv _{k+1} \ev _k ^T \zv _k.
\end{equation*}
Defining $\zv _k = [\zeta _1, ..., \zeta _k ] ^T$ and remembering that $\qv _{k+1} ^T \qv _{k+1} =1 $, we are led to
\begin{equation*}
    \| \bv - \Am ^T \uv _k \| _2  = \beta _ {k+1} \zeta _k.
\end{equation*}
We focus now on $ \zeta _k $. The matrix $\Lm _k ^T$ is a unit upper triangular (see \Cref{prop:KlowTri}), so
\begin{equation*}
    \zeta _k = \ev _k ^T \zv _k =
       \beta _1 \ev _k ^T \Lm _k ^{-T} \Bm _k ^{-T} \ev _1 = \beta _1 \ev _k ^T \Bm _k ^{-T} \ev _1.
\end{equation*}
Let us define
\begin{equation}
     \xv _k= \beta _1 \Bm _k ^{-T} \ev _1 = [\chi _1, ..., \chi _k] ^T. 
\end{equation}
Then, $\| \bv - \Am ^T \uv _k \| _2= \beta _ {k+1} \chi _k$.

At each iteration, the matrix  $\Bm _k ^T$ grows by a single row and column, with a new $\alpha _k$ and $\beta _k$, while its previous entries remain unchanged. Similarly, the vector $\beta _1 \ev _1$ only extends by an additional $0$. Then, $\xv _k$ is extended with one new entry per iteration. Under these circumstances, we can apply the same strategy used in \cite{arioli2013generalized} (and \cite{saad2003iterative}, for the DIOM version of \ac{FOM}) and compute only this new entry at each step following the recursion
\begin{equation*}
    \chi _1 =  \frac{\beta _1}{\alpha _1},  \qquad
    \chi _k = - \frac{\beta _k}{\alpha _k}  \chi _{k-1}.
\end{equation*}
\end{proof}

Although using the residual norm as stopping criterion for our iterative method is attractive for its simplicity and generality, it is important to also mention the following disadvantage. For challenging problems, there may be a gap between the residual norm and the error norm, making the former less reliable as stopping criterion.  We continue by providing an estimate of the energy norm for the error of the first solution vector $\uv$. This quantity can be used as alternative stopping criterion for \Cref{alg:nsGKB1}.

\begin{lemma}
\label{lem:errUnorm}
    The error for the primal variable at step $k$ 
    \begin{equation}
        \ev^{(k)} = \uv - \uv _k
    \end{equation}
    in the energy norm, is equal to
    \begin{equation}
        \normM{\ev^{(k)}} = \sqrt{ \xv _{n-k} ^T \Lm _{n-k} ^{-T} \xv _{n-k} }.
    \end{equation}
\end{lemma}
\begin{proof}
    First, considering \Cref{lem:nsGKBsol}, we have
    \begin{equation}
    \label{eq:LemPrfErr}
        \ev^{(k)} = \uv - \uv _k = \Vm _n \zv _n - \Vm _k \zv _k.
    \end{equation}
    The error is related to the quantities we have not yet found, from the iterations $k+1$ until $n$. For the remainder of this proof, we consider two indices: $k$, denoting quantities we already know, computed at iteration $k$ (\textbf{k}nown), and $n-k$, which refers to the \textit{last} $n-k$ (u\textbf{n-k}nown) entries from a quantity computed at iteration $n$. 
    With this notation, we have
    \begin{equation*}
        \xv _n = 
  [\chi _1, \dots, \chi _k, \chi _{k+1}, \dots, \chi _n ] ^T
    =
    \left[
	\begin{matrix}
		\xv _k \\
        \xv _{n-k} 
	\end{matrix}
	\right],
    \end{equation*}
    \begin{equation*}
        \Vm _n = [ \vvv _1, \dots, \vvv _k, \vvv _{k+1}, \dots, \vvv _n ] = [\Vm _k \ \Vm _{n-k} ],
    \end{equation*}
    \begin{equation*}
        \Vm _n ^T \Mm \Vm _n = \Lm _n =
        \left[
    	\begin{matrix}
    		\Lm _k & \mZ \\
            \Lm _{n-k,k} & \Lm _{n-k}
    	\end{matrix}
	    \right].
    \end{equation*}
    The vector $\zv$ not only gains a new entry per iteration, but the existing ones change as well. As such, we need to identify both the iteration and the entry with separate indices: 
    \begin{equation*}
        \zv _n = 
        [\zeta _1 ^{(n)}, \dots, \zeta _k ^{(n)}, \zeta _{k+1} ^{(n)}, \dots, \zeta _n ^{(n)}]  ^T=
        \left[
    	\begin{matrix}
    		\zv _k ^{(n)} \\
            \zv _{n-k} ^{(n)} 
    	\end{matrix}
	    \right]
     .
    \end{equation*}
    
We return to the proof, to the subtraction in \Cref{eq:LemPrfErr}. To do so, we write $\zv _n = \Lm _n ^{-T} \xv _n$ in more detail, using the notation introduced above. Since $\Lm _k$ and $\Lm _{n-k}$ are invertible, we obtain
 \begin{equation*}
     \left[
    \begin{matrix}
        \zv _k ^{(n)} \\
        \zv _{n-k} ^{(n)} 
    \end{matrix}
    \right]=
     \left[
    \begin{matrix}
        \Lm _k & \mZ \\
        \Lm _{n-k,k} & \Lm _{n-k}
    \end{matrix}
    \right] ^{-T} 
     \left[
	\begin{matrix}
		\xv _k \\
        \xv _{n-k} 
	\end{matrix}
	\right]=
    \left[
    \begin{matrix}
        \Lm _k ^{-T}  & \Ym \\
        \mZ  & \Lm _{n-k} ^{-T} 
    \end{matrix}
    \right] 
     \left[
	\begin{matrix}
		\xv _k \\
        \xv _{n-k} 
	\end{matrix}
	\right],
 \end{equation*}
 with $\Ym = - \Lm _k ^{-T} \Lm _{n-k,k} ^T \Lm _{n-k} ^{-T}$. 
 Then,
 \begin{equation}
 \label{eq:link_zn_zk}
    \zv _k ^{(n)} = \Lm _k ^{-T} \xv _k + \Ym \xv _{n-k} \ \text{ and } \ \zv _{n-k} ^{(n)} = \Lm _{n-k} ^{-T} \xv _{n-k}.
 \end{equation}
Note that $\zv _k = \Lm _k ^{-T} \xv _k$, which, considering the above, allows us to connect $\zv _k$ and the first part of $\zv _n$ through $\zv _k ^{(n)} = \zv _k + \Ym \xv _{n-k}$. Then, it follows for \Cref{eq:LemPrfErr} 
\begin{equation*}
    \begin{split}
        \ev^{(k)} &= \Vm _n \zv _n - \Vm _k \zv _k = 
        \left[ \begin{matrix}
            \Vm _k & \Vm _{n-k} 
        \end{matrix} \right] 
        \left[
        \begin{matrix}
            \zv _k ^{(n)} \\
            \zv _{n-k} ^{(n)} 
        \end{matrix}
        \right] 
        - \Vm _k  \zv _k \\
        &= 
        \Vm _k (\zv _k ^{(n)} - \zv _k) + \Vm _{n-k} \zv _{n-k} ^{(n)} \\
        &= \Vm _k \Ym \xv _{n-k} + \Vm _{n-k} \Lm _{n-k} ^{-T} \xv _{n-k} \\        
        &= \Vm _n 
        \left[
        \begin{matrix}
            \Ym \\
            \Lm _{n-k} ^{-T}
        \end{matrix}
        \right] \xv _{n-k}.
    \end{split}
\end{equation*}
We continue with the (squared) norm of the error vector.
\begin{equation*}
\begin{split}
    \normM{\ev^{(k)} } ^2 &= \ev^{(k) T} \Mm \ev^{(k)} \\    
    &= \xv _{n-k} ^T 
    \left[ \begin{matrix}
        \Ym ^T & \Lm _{n-k} ^{-1} 
    \end{matrix} \right] 
     \left[
    \begin{matrix}
        \Lm _k & \mZ \\
        \Lm _{n-k,k} & \Lm _{n-k}
    \end{matrix}
    \right]
     \left[
    \begin{matrix}
        \Ym \\
        \Lm _{n-k} ^{-T}
    \end{matrix}
    \right]
    \xv _{n-k}
\end{split}
\end{equation*}

Using the definition of $\Ym$ above, the matrix product simplifies to
\begin{flalign*}
    & \left[ \begin{matrix}
        \Ym ^T & \Lm _{n-k} ^{-1} 
    \end{matrix} \right] 
     \left[
    \begin{matrix}
        \Lm _k & \mZ \\
        \Lm _{n-k,k} & \Lm _{n-k}
    \end{matrix}
    \right]
     \left[
    \begin{matrix}
        \Ym \\
        \Lm _{n-k} ^{-T}
    \end{matrix}
    \right] = &
\end{flalign*}

 \begin{flalign*}   
     & -  \Lm _{n-k} ^{-1} \Lm _{n-k,k} \underbrace{\Lm _k ^{-1} \Lm _k} _{\Id}  \Ym + \Lm _{n-k} ^{-1} \Lm _{n-k,k} \Ym + \underbrace{\Lm _{n-k} ^{-1} \Lm _{n-k} } _{\Id} \Lm _{n-k} ^{-T} = \Lm _{n-k} ^{-T}. &
\end{flalign*}
Finally, the squared energy norm of the error is 
\begin{equation}
    \normM{\ev^{(k)}} ^2= \xv _{n-k} ^T \Lm _{n-k} ^{-T} \xv _{n-k}.
\end{equation}
\end{proof}

It is also possible to stop the iterative process using an estimate for the energy norm of the error associated with the primal variable $ \uv $ similar to the one described in \cite{arioli2013generalized}, using the last $d$ entries of the vectors $\zv _k$ and $\xv _k$. In this case, we need to compute $d$ steps in the backwards substitution with $\Lm _k ^T$ at each iteration.  Given $\zv _k = [\zeta _1, ..., \zeta _k ] ^T= \Lm _k ^ {-T}[\chi _1, ..., \chi _k ] ^T$, with $ k \geq d $, the error norm estimate corresponding to step $k - d $ is
\begin{equation}
\label{eq:errLowBnd}
    \xi _{k-d} ^2 = \sum _{i=k-d+1} ^{k} \chi _i \zeta _i.
\end{equation}
We can adapt this, in order to use it as a relative stopping criterion
\begin{equation}
\label{eq:errBnd}
    \frac {    \xi _{k-d} ^2 }  { \sum _{i=1} ^{k} \chi _i \zeta _i }.
\end{equation}

The delay parameter $d$ acts as an additional safety measure, leading to a few extra iterations. In practical experiments, $d$ is heuristically chosen to be e.g. 5, as in \cite{arioli2013generalized,KrSoArTaRu2021}.

One should also bear in mind that since \ac{FOM} does not feature a minimization property, neither does our solver, in light of the equivalence described in \Cref{sec:nsGKB}. This implies  that we cannot guarantee a decrease in the residual/error norm at each iteration of the algorithm.

\import{./}{new_tests.tex}

\section{Conclusions}
\label{sec:conc}

In this work we have developed nsCRAIG, a solver for nonsymmetric saddle point systems. 
Its theoretical foundation is inspired by a known equivalence for the symmetric case, between generalized CRAIG \cite{arioli2013generalized} and CG. 
We extended this such that nsCRAIG is equivalent to FOM instead. 
Aside from the theoretical point of view, we have illustrated this relationship also in a numerical setting via our experiments.

Along with our solver's description, proof of convergence and algorithm, we also provided stopping criteria. 
One is an inexpensive way to compute the residual norm for the second equation of the saddle point system. 
The other choice is an estimate of the error for the primal variable in an energy norm. 

We compared our solver with \ac{GMRES}, a popular choice for tackling nonsymmetric and indefinite problems. 
The findings confirmed our expectations concerning the lower memory requirements of our algorithm. 
The difference stems from a key feature: nsCRAIG generates and stores a right basis with shorter vectors compared to \ac{GMRES}. 
The left basis from nsCRAIG (with long vectors) does not need to be stored, significantly decreasing memory costs. This advantage makes our method an attractive choice, especially in cases where \ac{GMRES} would need to be restarted often, with the known negative impact on convergence.

We have also found that after suitable preconditioning, our algorithm is at least twice as fast compared to \ac{GMRES} without restarts. This is due to the latter's convergence behavior, where often only every other iteration significantly progresses towards convergence. A similar pattern was already known for the symmetric case, between generalized CRAIG and MINRES \cite{arioli2013generalized}.

In this work, we have only considered exact matrix vector products of type $\Mm ^{-1} \vvv$ involving the inverse of the leading block of the saddle point system. A more general and practically motivated alternative is to consider an inexact approach by making use of an iterative solver for this inner problem. For the symmetric generalized variant of CRAIG, such a strategy has been explored in \cite{darrigrand2022inexact} and yielded promising results. A similar study concerning the method we presented in this paper could constitute an interesting direction for future developments.

\bibliographystyle{plain}
\bibliography{references}

\end{document}

%% file: new_tests.tex
\section{Numerical experiments}
\label{sec:tests}

Our application of choice is from the field of Computational Fluid Dynamics, in the form of Navier-Stokes flow problems. 
To generate the linear systems for our tests, we make use of the Incompressible Flow \& Iterative Solver Software \footnote{http://www.cs.umd.edu/~elman/ifiss3.6/index.html} (IFISS) package (see also \cite{elman2014ifiss}).
The particular problems under consideration are those also used in \cite{elman2008least}, and described in more detail in \cite{elman2014finite}. 
We briefly summarize them here. 
In each case, we generate the nonlinear Navier-Stokes problem with IFISS, which is given by 
 \begin{align}
  - \nu \nabla ^2 \Vec{u} + \Vec{u}   \cdot \nabla \Vec{u} + \nabla p &= \Vec{f} \label{eq:contNavSto}\\
 	\nabla \cdot  \Vec{u} &=0,
 \end{align}
with the kinematic viscosity $\nu$. 
To deal with the nonlinearity of the Navier--Stokes equation in the convection term, \textit{Picard's iteration} is used to obtain the following linearized  equations 
\begin{align}
    {-\nu\Delta\vec{u}^{(k)}}+{(\vec{u}^{(k-1)}\cdot\nabla)\vec{u}^{(k)}}+{\nabla p^{(k)}}&={\vec{f}}&\text{in }\Omega,\label{oseen:eq1}\\
    {\nabla\cdot\vec{u}^{(k)}}&=0&\text{in }\Omega,\label{oseen:eq2}
\end{align}
for each iteration $k$, starting from an arbitrary initial guess $(\vec{u}^{(0)},p^{(0)})$. For the linearized system we consider a Q2-Q1 Finite Element discretization, leading to a saddle point system where the leading block is nonsymmetric and positive definite (see \cite{konshin2015ilu}).

Next, we briefly describe the two test problems. Picard's method iterates until the nonlinear residual norm reaches $10^{-5}$. The next linear system used for the correction step represents the input for the linear solvers we consider in \Cref{subsec:numerComp}. 

\paragraph{Test case 1: Flow over a backward facing step}
This case represents a flow with a parabolic inflow velocity profile passing through a domain $\Omega = ( (-1, 5 ) \times (-1, 1 ) ) \setminus ( (-1, 0 ] \times (-1, 0] )$. The boundary conditions are

\begin{equation*}
    \begin{dcases*}
u_x = 4y(1-y), \quad u_y =0  \qquad  
& at the inflow $ \Gamma_{in}= \left\{-1\right\} \times \left[0, 1 \right]$,\\
\nu \frac{\partial u_x}{\partial x} - p = 0, \quad \frac{\partial u_y}{\partial x}  = 0 \qquad 
& at the outflow $ \Gamma_{out}= \left\{5 \right\} \times \left[-1, 1 \right]$, \\ 
\text{no-slip } \qquad & on the horizontal walls.
    \end{dcases*}   
\end{equation*}

\paragraph{Test case 2: Driven cavity flow}
The domain for this problem is $\Omega =  (-1, 1 ) \times (-1, 1 ) $, with the following boundary conditions

\begin{equation*}
    \begin{dcases*}
u_x= 1-x^4, \quad u_y=0 \qquad
& on the wall $ \Gamma_{top}=   \left[-1, 1 \right] \times \left\{-1\right\} $,\\
\text{no-slip } \qquad & on the bottom and vertical walls.
    \end{dcases*}   
\end{equation*}
This represents a model where the cavity lid is moving according to the given regularized condition, driving the enclosed flow.

\subsection{Comparison with \ac{FOM} and \ac{GMRES}}
\label{subsec:numerComp}

Saddle-point systems are often solved with the \ac{MINRES} method when the leading block is symmetric and with \ac{GMRES} when this condition is not fulfilled. These approaches treat the system matrix as a whole, i.e., in an all-at-once manner. As such, \ac{GMRES} stores vectors of length $m+n$ in memory. 
In contrast, using the approach we presented in \Cref{sec:nsGKB}, only the right vectors of size $n$ need to be stored. Depending on the nature of the problem and the associated discretization, $n$ may be (much) smaller than $m$. As consequence, for a given amount of memory, we could store more vectors than in the case of \ac{GMRES}. We might hence limit the negative effects on the number of iterations observed in strategies like restarting.

In this section, we compare our proposed approach to \ac{GMRES} by measuring the necessary number of iterations to reach convergence, defined as reducing the relative residual norm below $10^{-3}$. Additionally, we compare the memory usage of both algorithms. In a first scenario, we impose no memory limitation. Following the considerations above, we expect \ac{GMRES} to be more expensive in this regard. In a second scenario, we limit the amount of memory \ac{GMRES} can use such that it does not exceed that used by nsCRAIG in \Cref{alg:nsGKB1}. Once GMRES reached the available amount of memory usage, it will be restarted every $k _{\textit{max}}$ iterations. This value is computed as
\begin{equation}\label{eq:restart_kmax}
    k _{\textit{max}} = \floor*{ \frac{\textit{iter} \cdot n}{m+n}  } .
\end{equation}

Our choice follows from an estimate of the memory used by \Cref{alg:nsGKB1} after $\textit{iter}$
steps:
\begin{equation}
\label{eq:gkb_memory}
    m + n (\textit{iter}+1),
\end{equation}
where $m$ is the length of a vector holding either $\wv, \vvv$ or $\uv$, then $n$ for a vector holding $\bv, \gvv, \yv$ or $\pv$ and finally $n * \textit{iter}$ for all the $\qv$ vectors. We consider the approximate memory usage of a textbook implementation of restarted \ac{GMRES} (see, e.g., \cite{saad2003iterative}) to be
\begin{equation}
\label{eq:gmres_memory}
    (m + n) (k _{\textit{max}} + 1).
\end{equation}
By setting \Cref{eq:gmres_memory} equal to \Cref{eq:gkb_memory} and solving for $k _{\textit{max}} $, we arrive at \Cref{eq:restart_kmax}.
For the following comparisons, we take \Cref{alg:nsGKB1} as our representative implementation, since it matches a typical use case. It computes the explicit approximate solution only at the end, after the relative residual norm has reached the given threshold. 

Furthermore, in \Cref{sec:nsGKB}, we have based our developments on a theoretical equivalence between \ac{FOM} on the Schur complement and our algorithm. 
In this section, we test whether the  equivalence also holds numerically.

Since we compare the number of iterations performed by each solver, it is important that all of them have a comparable cost. In this regard, the most expensive operation in FOM and our algorithm is the need to apply $\Mm ^{-1}$ to a vector. 
To bring \ac{GMRES} to the same level, we consider it as applied to the 
right-preconditioned problem
\begin{equation}\label{eq:prec_sys}
	\left[
	\begin{matrix}
		\Mm & \Am \\
		\Am ^T & \mZ  
	\end{matrix}
	\right]
    \left[
	\begin{matrix}
		\Mm ^{-1} & \mZ  \\
		\mZ  & \Id 
	\end{matrix}
	\right]
	\left[
	\begin{matrix}
		\tilde{\uv}  \\
		\tilde{\pv} 
	\end{matrix}
	\right]
	=
	\left[
	\begin{matrix}
		\mZ \\
		\bv 
	\end{matrix}
	\right],
\end{equation}
\begin{equation}\label{eq:unprec_sol}
    \left[
	\begin{matrix}
		\Mm  & \mZ  \\
		\mZ  & \Id 
	\end{matrix}
	\right]
	\left[
	\begin{matrix}
		\uv \\
		\pv
	\end{matrix}
	\right]
	=
	\left[
	\begin{matrix}
		\tilde{\uv}  \\
		\tilde{\pv} 
	\end{matrix}
	\right] .
\end{equation}
The $\Mm ^{-1}$ in \Cref{eq:prec_sys} is applied every iteration, while the system in \Cref{eq:unprec_sol} is solved only once, at the end, to find the solution of the unpreconditioned \Cref{eq:ini_sys}.

We plot the convergence history of the linear solvers in \Cref{fig:NSStepGrid5Visc1Ov100_1e3,fig:NSCavityGrid5Visc1Ov200_1e3}, for the step and cavity problems respectively. 
It is visible how \ac{FOM} and our \Cref{alg:nsGKB1} behave identically in both cases and are significantly faster than unrestarted right-preconditioned \ac{GMRES}. 
This is further emphasized by comparing against a resatrted version of \ac{GMRES}.  
Here the restart 
parameter has been computed according to \Cref{eq:restart_kmax}. 
If constrained to use the same amount of memory, \ac{GMRES} needs about five times more iterations than our algorithm to reach convergence. 
In terms of memory usage, GMRES without restart needs 17 time more than our solver for the step problem and 24 times more for the cavity problem.  

\import{images/pgf/}{NSStepGrid5Visc1Ov100_1e3.tex}

It is interesting to note that \ac{GMRES} without restarts has phases where only every other iteration significantly contributes to the objective of reducing the residual norm (see \Cref{fig:zoomNSStep}). 
A similar behavior has already been noted in the symmetric case, when the generalized Golub-Kahan is compared with \ac{MINRES} \cite{arioli2013generalized}. 
The explanation is related to the particular choice of block preconditioner used (see \Cref{eq:prec_sys}), which leads to a matrix with a symmetric spectrum. 
Such a spectrum has an impact on the convergence of the solver, which reduces the residual only every other step (see \cite{fischer1998minimum}).

\import{images/pgf/}{NSCavityGrid5Visc1Ov200_1e3.tex}

\import{images/pgf/}{zoomNSStep.tex}

%% file: images/pgf/NSStepGrid5Visc1Ov100_1e3.tex
\begin{figure}[H]
    \centering
    \begin{tikzpicture}
    \begin{axis}[
	table/col sep=tab,
    ymode=log,    
    legend pos=north east,
    ytick={1,1e-1,1e-2,1e-3},
    width=\FigWid \textwidth,
	height=\FigHei \textwidth,
    xlabel={Iterations},
    ylabel={Relative residual norm}
    ]
	
\addplot+[] table [x=Nonsym. GKBX, y=Nonsym. GKBY]{images/csv/NSStepGrid5Visc1Ov100_1e3.csv}; 
\label{fig:item:NSStepGrid5Visc1Ov100_1e3_Nonsym. GKB} 
\addlegendentry{Nonsym. GKB} 

\addplot+[] table [x=FOM Schur Comp.X, y=FOM Schur Comp.Y]{images/csv/NSStepGrid5Visc1Ov100_1e3.csv}; 
\label{fig:item:NSStepGrid5Visc1Ov100_1e3_FOM Schur Comp.} 
\addlegendentry{FOM Schur Comp.} 

\addplot+[] table [x=GMRESX, y=GMRESY]{images/csv/NSStepGrid5Visc1Ov100_1e3.csv}; 
\label{fig:item:NSStepGrid5Visc1Ov100_1e3_GMRES} 
\addlegendentry{GMRES} 

\addplot+[] table [x=GMRES(15)X, y=GMRES(15)Y]{images/csv/NSStepGrid5Visc1Ov100_1e3.csv}; 
\label{fig:item:NSStepGrid5Visc1Ov100_1e3_GMRES(15)} 
\addlegendentry{GMRES(15)} 

\end{axis} 
\end{tikzpicture} 
\caption{Convergence history of the listed solvers for the step domain test case.} 
\label{fig:NSStepGrid5Visc1Ov100_1e3} 
\end{figure}
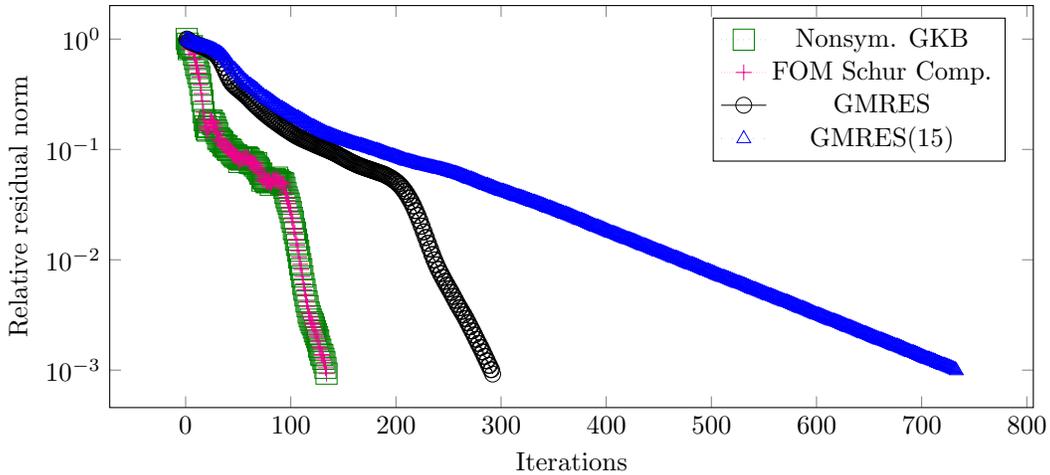

%% file: images/pgf/NSCavityGrid5Visc1Ov200_1e3.tex
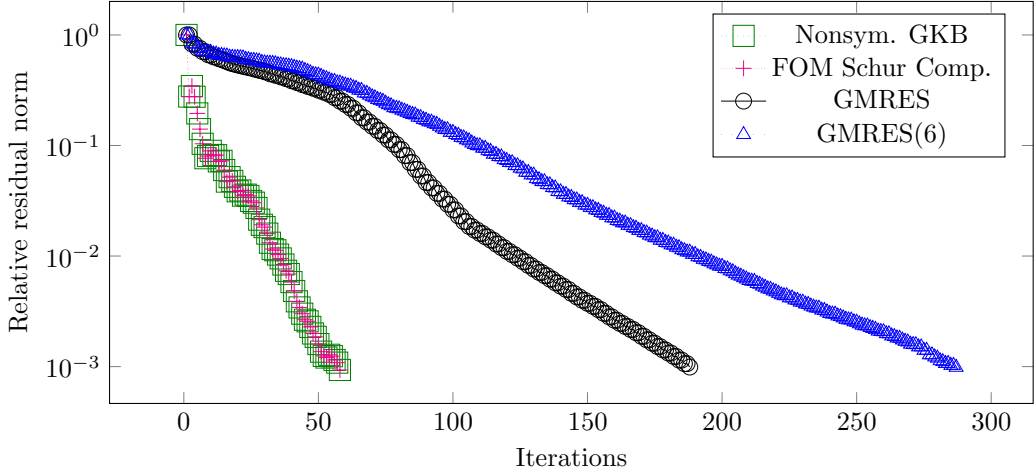
\begin{figure}[H]
    \centering
    \begin{tikzpicture}
    \begin{axis}[
	table/col sep=tab,
    ymode=log,    
    legend pos=north east,
    ytick={1,1e-1,1e-2,1e-3},
    width=\FigWid \textwidth,
	height=\FigHei \textwidth,
    xlabel={Iterations},
    ylabel={Relative residual norm}
    ]
	
\addplot+[] table [x=Nonsym. GKBX, y=Nonsym. GKBY]{images/csv/NSCavityGrid5Visc1Ov200_1e3.csv}; 
\label{fig:item:NSCavityGrid5Visc1Ov200_1e3_Nonsym. GKB} 
\addlegendentry{Nonsym. GKB} 

\addplot+[] table [x=FOM Schur Comp.X, y=FOM Schur Comp.Y]{images/csv/NSCavityGrid5Visc1Ov200_1e3.csv}; 
\label{fig:item:NSCavityGrid5Visc1Ov200_1e3_FOM Schur Comp.} 
\addlegendentry{FOM Schur Comp.} 

\addplot+[] table [x=GMRESX, y=GMRESY]{images/csv/NSCavityGrid5Visc1Ov200_1e3.csv}; 
\label{fig:item:NSCavityGrid5Visc1Ov200_1e3_GMRES} 
\addlegendentry{GMRES} 

\addplot+[] table [x=GMRES(6)X, y=GMRES(6)Y]{images/csv/NSCavityGrid5Visc1Ov200_1e3.csv}; 
\label{fig:item:NSCavityGrid5Visc1Ov200_1e3_GMRES(6)} 
\addlegendentry{GMRES(6)} 

\end{axis} 
\end{tikzpicture} 
\caption{Convergence history of the listed solvers for the cavity domain test case.} 
\label{fig:NSCavityGrid5Visc1Ov200_1e3} 
\end{figure}

%% file: images/pgf/zoomNSStep.tex
\begin{figure}[H]
    \centering
    \begin{tikzpicture}
    \begin{axis}[
	table/col sep=tab,
    ymode=log,
    restrict x to domain=270:300,
    legend pos=north east,
    xtick={270,280,290},
    width=\FigWid \textwidth,
	height=\FigHei \textwidth,
    xlabel={Iterations},
    ylabel={Relative residual norm}
    ]

\addplot+[color=black, densely dotted, mark=o, line width=1.0pt, mark size=3pt] table [x=GMRESX, y=GMRESY]{images/csv/NSStepGrid5Visc1Ov100_1e3.csv}; 
\label{fig:item:zoomNSStep_GMRES} 
\addlegendentry{GMRES} 

\end{axis} 
\end{tikzpicture} 
\caption{Last 30 \ac{GMRES} iterations for the step domain test case. Zoom-in of \Cref{fig:NSStepGrid5Visc1Ov100_1e3}.} 
\label{fig:zoomNSStep} 
\end{figure}
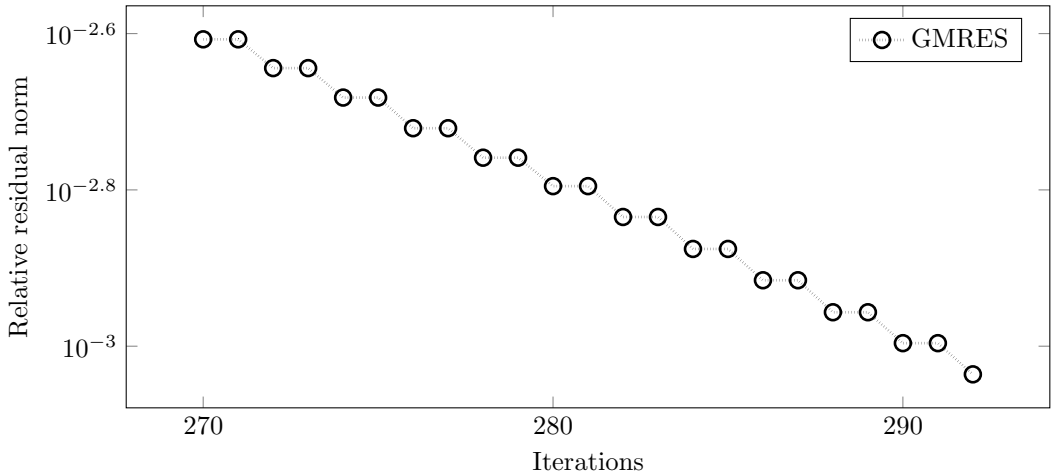